\documentclass[11pt,oneside, english]{article}
\usepackage{dsfont}
\usepackage{amsmath, amssymb, amsthm, verbatim,enumerate, bbm}
\usepackage{enumerate}
\usepackage[utf8x]{inputenc}
\usepackage[T1]{fontenc}
\usepackage[english]{babel}
\usepackage{mathtools}
\usepackage[
letterpaper,
top=0.8in,
bottom=0.8in,
left=1in,
right=1in]{geometry}
\usepackage{bm} %
\usepackage[colorlinks=true, allcolors=blue]{hyperref}

\usepackage{graphicx}

\graphicspath{{./figures/}}
\usepackage{epsfig}
\usepackage{fullpage}
\usepackage{mathrsfs}
\usepackage{enumerate}
\usepackage{verbatim}
\usepackage{mhequ}
\usepackage{algorithm}
\usepackage{algorithmicx}
\usepackage{algpseudocode}
\usepackage[nameinlink,capitalise,noabbrev]{cleveref}

\usepackage{showlabels}
\numberwithin{equation}{section}

\def \be{\begin{equs}}
	\def \ee{\end{equs}}

\def \be{\begin{equs}}
\def \ee{\end{equs}}

\newtheorem*{namedtheorem}{\theoremname}
\newcommand{\theoremname}{testing}

\newtheorem*{question*}{Question}

\theoremstyle{definition}

\theoremstyle{plain}

\title{Elephant Random Walk on Triangular Lattice}

\author{
    Rohit Chaudhuri\thanks{Indian Statistical Institute, Kolkata}
    }

\date{}

\begin{document}

\pagestyle{empty}

\maketitle
\section{Abstract}

In this report, we introduce the elephant random walk on the triangular lattice over $R^2$ incorporating directions by extending the model developed in \cite{baur2016elephant}. We study the behaviour of the walk by finding the appropriate scaling limits.\medskip

\section{The walk}

Let $\{X_n\}_{n\in\mathbb{N} }$  denote the step size random variables for the walk. Now, $X_1 \sim Uniform \{\pm 1,\pm \omega, \pm \omega^2 \}$ where $\omega$ is the standard notation for cube root of unity. For each $n\in\mathbb{N}$ we do the following: \newline Choose $T_n \sim Uniform \{1, \dots ,n \}$. Then,

\begin{center}
$X_{n+1}= X_{T_n} $, w.p $p/2$

\hspace{1cm}$= -X_{T_n}$,  w.p $p/2$

\hspace{1cm}$= \xi_{n+1} $, w.p $ 1-p $

\end{center}

where $\{\xi_n\}_{n\in\mathbb{N} } \stackrel{i.i.d}{\sim} Uniform \{\pm 1,\pm \omega, \pm \omega^2 \} $. \newline Now, define $S_n = \displaystyle\sum_{i=1}^{n} X_i$ for each $n\in\mathbb{N}$ (where $S_{0}=0$). Then, $\{S_n\}_{n\in\mathbb{N} }$ denotes the required walk. 

\section{Results}

Suppose, $\{S_n\}_{n\in\mathbb{N} }$ is the walk as introduced above. Then,

\begin{itemize}
	\item $\frac{S_{n}}{n} \stackrel{n \to \infty}{\to} 0$ a.e.
	\item $\frac{S_{n}}{\sqrt{n}} \stackrel{n \to \infty}{\to} Y$ in distribution. 
\end{itemize}
where the random variable $Y$ has a gaussian distribution independent of the parameter $p$.
\section{Proofs}

For each $n\in\mathbb{N} $, let $U_{n}$ be the indicator of whether random renovation (i.e. $\xi_{n+1}$) is required for constructing $X_{n+1}$. Then, $\{U_n\}_{n\in\mathbb{N} } \stackrel{i.i.d}{\sim} Ber(1-p)$.\\ Now, define $\{Y_n\}_{n\in\mathbb{N} }$ and $\{Z_n\}_{n\in\mathbb{N} }$  such that:

\begin{center}
	$Y_{n}= 1 $, if $X_{n} \in\{ 1, \omega, \omega^2\}$
	
	\hspace{1.5cm}$= -1$,  if $X_{n} \in\{ -1,-\omega,-\omega^2 \}$
		
\end{center}

and $Z_{n}=X_{n}/{Y_n}$ for each $n\in\mathbb{N}$. Clearly, $\{Z_n\}_{n\in\mathbb{N} }\in\{1,\omega,\omega^2\}$.\\[1cm] \textbf{Claim 1:} $\{Z_n\}_{n\in\mathbb{N} }\amalg \{Y_n\}_{n\in\mathbb{N} }$ and $\{Y_n\}_{n\in\mathbb{N} }\stackrel{i.i.d}{\sim}Ber(1/2)$ .\\[0.5cm] \textit{Proof:}  Choose $\{y_m\}_{m\in\{1,2,\dots,n+1\} }\in\{+1,-1\}$ , $\{z_m\}_{m\in\{1,2,\dots,n\} }\in\{ 1, \omega, \omega^2\}$ and $t\in\{1,\dots,n\}$.\\[0.2cm]$P[Y_{n+1}=y_{n+1}|\{Y_m\}_{m\in\{1,2,\dots,n\}}=\{y_m\}_{m\in\{1,2,\dots,n\} },\{Z_m\}_{m\in\{1,2,\dots,n\} }=\{z_m\}_{m\in\{1,2,\dots,n\} },T_{n}=t]$ \\[0.2cm] $=P[Y_{n+1}=y_{n+1}|Y_{t}=y_{t}, Z_{t}=z_{t}]$ \\[0.2cm] $=P[Y_{n+1}=y_{n+1}, U_{n}=0|Y_{t}=y_{t}, Z_{t}=z_{t}]+P[Y_{n+1}=y_{n+1}, U_{n}=1|Y_{t}=y_{t}, Z_{t}=z_{t}]$ \\[0.2cm] $=\frac{p}{2}+\frac{1-p}{2}$ \\[0.2cm] $=\frac{1}{2}$. $\Box$\\[1cm]\textbf{Claim 2}: Let $C_{n}^{i}= \displaystyle\sum_{m=1}^{n} I_{Z_{m}=\omega^i}$ for $i=1,2,3$, then $\frac{C_{n}^{i}}{n}\xrightarrow[n\to \infty]{a.s}\frac{1}{3}$ for $i=1,2,3$. \\[0.2cm]\textit{Proof:} Clearly we have:
\begin{center}
$Z_{n+1}= Z_{T_n} $, w.p $p$

\hspace{2cm}$=\mid\xi_{n+1}\mid$, w.p $ 1-p $
\end{center}

where $\mid.\mid$ ignores the $\pm$ sign.\\[1cm] Hence $\{Z_n\}_{n\in\mathbb{N} }$ constitutes the step sizes of ERW over $\{1,\omega,\omega^2\}$. Let $C_n=(C_{n}^{1},C_{n}^{2},C_{n}^{3})$ denote the count vector. We will try to show that $\{C_n\}_{n\in\mathbb{N}}$ evoles like the color count of urn model involving three colors (representing $\omega$, $\omega^2, 1$). Before that, note that $\sigma(C_{1},C_{2},\dots,C_{n})=\sigma(Z_{1},Z_{2},\dots,Z_{n})\hspace{0.2cm}\forall n\in \mathbb{N}$. Let, $e_1=(1,0,0), e_2=(0,1,0), e_3=(0,0,1)$ and $D_{n}=C_{n+1}-C_{n}\hspace{0.2cm} \forall n\in\mathbb{N}$. So,\\[1cm]$P[D_n=e_1|C_1,\dots,C_n]$ \\[0.2cm] $=P[D_n=e_1|Z_1,\dots,Z_n]$ \\[0.2cm] $=\displaystyle\sum_{i=1}^{n} P[D_n=e_1,T_n=i|Z_1,\dots,Z_n]$\\[1cm] $=\frac{1}{n}\displaystyle\sum_{i=1}^{n} P[D_n=e_1|Z_1,\dots,Z_n,T_n=i]$ \\[0.2cm]$=\frac{1}{n}\displaystyle\sum_{i=1}^{n} P[D_n=e_1|Z_i,T_n=i]$ \\[0.5cm]Choose any $i\in\{1,\dots n\}$, suppose $Z_{i}=\omega^j$, where $j\neq1$. Then,\\[0.2cm]$P[D_{n}=e_1|Z_{i}=\omega^j,T_n=i]$ \\[0.2cm]$=P[\xi_{n+1}=\pm\omega,U_{n}=0|Z_{i}=\omega^j,T_n=i]$\\[0.2cm]$=P[\mid\xi_{n+1}\mid=\omega,U_{n}=0]$\\[0.2cm]$=\frac{1-p}{3}$. \\[0.5cm] Now, suppose $Z_i=\omega$. Then, \\[0.2cm]$P[D_{n}=e_1|Z_{i}=\omega,T_n=i]$\\[0.2cm]$=P[U_{n}=0|Z_{i}=\omega,T_n=i]+P[\xi_{n+1}=\pm\omega,U_{n}=1|Z_{i}=\omega,T_n=i]$\\[0.2cm]$=p+\frac{1-p}{3}$.\\[0.5cm]Hence, combining all these and continuing the previous calculations we get:\\[0.2cm]$P[D_n=e_1|C_1,\dots,C_n]$\\[0.2cm]$=\frac{1}{n}\displaystyle[\displaystyle\sum_{i:Z_i=\omega} P[D_n=e_1|Z_i,T_n=i]+\displaystyle\sum_{i:Z_i\neq \omega} P[D_n=e_1|Z_i,T_n=i]\displaystyle]$\\[0.2cm]$=\frac{1}{n}\displaystyle[\displaystyle\sum_{i:Z_i=\omega} (p+\frac{1-p}{3})+\displaystyle\sum_{i:Z_i=\omega^2} (\frac{1-p}{3})+\displaystyle\sum_{i:Z_i=\omega^3} (\frac{1-p}{3})\displaystyle]$\\[0.2cm]$=(p+\frac{1-p}{3})\frac{C_{n}^1}{n}+(\frac{1-p}{3})\frac{C_{n}^2}{n}+(\frac{1-p}{3})\frac{C_{n}^3}{n}$.\\[1cm]Similarly, we can compute $P[D_n=e_2|C_1,\dots,C_n]$, $P[D_n=e_3|C_1,\dots,C_n]$ and from these transition probabilities we can say that $\{C_n\}_{n\in\mathbb{N}}$ evoles like the color count of urn model involving three colors (representing $\omega$, $\omega^2, 1$) with i.i.d random replacement matrices $\{R_n\}_{n\in\mathbb{N}}$ sattisfying:

\[
E(R_n)=
\begin{pmatrix}
p+\frac{1-p}{3} & \frac{1-p}{3} & \frac{1-p}{3} \\
\frac{1-p}{3} & p+\frac{1-p}{3} & \frac{1-p}{3} \\
\frac{1-p}{3} & \frac{1-p}{3} & p+\frac{1-p}{3}

\end{pmatrix}
,\forall n\in\mathbb{N}\] 
 
Since, the above matrices are doubly stochastic, owing to the theorem of Bai and Hu which is very standard, we get the desired result. $\Box$\\[1cm]

Now for each $i\in\{1,2,3\}$, we set $\tau_{0}^i=0$ and for each $n\in\mathbb{N}$ we define $\tau_{n}^i=inf\{m>\tau_{n-1}^i \mid Z_m=\omega^i\}$. From \textbf{claim 1}, we get that $\{Y_{\tau_{n}^1}\}_{n\in\mathbb{N}}$, $\{Y_{\tau_{n}^2}\}_{n\in\mathbb{N}}$, $\{Y_{\tau_{n}^3}\}_{n\in\mathbb{N}}$ are three independent i.i.d Ber(1/2) sequences which are independent of the sequence $\{C_n\}_{n\in\mathbb{N}}$. Also, note that $S_n = \displaystyle\sum_{i=1}^{3} \omega^i \displaystyle\sum_{m=1}^{C_n^i} Y_{\tau_{m}^i}$ from which the desired conclusion follows.

\phantomsection
  \addcontentsline{toc}{section}{References}
  \bibliographystyle{amsalpha}
  \bibliography{biblio.bib}
\end{document}